\documentclass{amsart}
\usepackage[english]{babel}
\usepackage[latin1]{inputenc}
 \usepackage[all]{xy}
\usepackage{amsmath,amsfonts,amssymb,amsthm,amscd,array,stmaryrd,mathrsfs, mathdots}
\usepackage{pstricks}
\usepackage{graphics, graphicx}

\theoremstyle{plain}
\newtheorem{thm}{Theorem}%[section]
\newtheorem{lem}{Lemma}[section]

\newtheorem{prop}[lem]{Proposition}
\newtheorem{conj}[]{Conjecture}

\theoremstyle{definition}
\newtheorem{defn}[lem]{Definition}
\newtheorem{rem}[lem]{Remark}
\newtheorem{ex}[lem]{Example}

\newcommand{\R}{\mathbb{R}}
\newcommand{\Z}{\mathbb{Z}}

\newcommand{\Q}{\mathbb{Q}}

\newcommand{\F}{\mathcal{F}}

\newcommand{\Pc}{\mathcal{P}}

\newcommand{\Rc}{\mathcal{R}}
\newcommand{\Sc}{\mathcal{S}}

\newcommand{\cX}{\mathcal{X}}

\newcommand{\SL}{\mathrm{SL}}
\newcommand{\PSL}{\mathrm{PSL}}

\begin{document}

\title[Quantum numbers and $q$-friezes]{Quantum numbers and $q$-deformed Conway-Coxeter friezes}

\author{Sophie Morier-Genoud}
\author{Valentin Ovsienko}

\address{Sophie Morier-Genoud,
Sorbonne Universit\'e and Universit\'e de Paris, CNRS, IMJ-PRG, F-75006 Paris, France
}

\address{
Valentin Ovsienko,
Centre national de la recherche scientifique,
Laboratoire de Math\'ematiques
U.F.R. Sciences Exactes et Naturelles
Moulin de la Housse - BP 1039
51687 Reims cedex 2,
France}
\email{sophie.morier-genoud@imj-prg.fr, valentin.ovsienko@univ-reims.fr}

\maketitle

\hfill{ \it To the memory of John Conway}

\thispagestyle{empty}

\bigskip

The notion of $q$-deformed real numbers was recently introduced in~\cite{SVRat,SVRe}. The first steps of this theory were largely influenced by John
Conway. In the fall of 2013, we had a chance to spend a week in his company.  The occasion was a conference that we organized in Luminy, at which
Conway gave two wonderful talks. We cannot forget the long conversations we had with John there; between desperate attempts to teach us his famous
doomsday algorithm, he told us much about his work with Coxeter. Our understanding of Conway and Coxeter's theorem has transformed into a certain
combinatorial viewpoint on continued fractions~\cite{MGO}, that eventually led us to $q$-deformations.

One goal of this article is to explain the origins of $q$-numbers and sketch their main properties. We will also introduce $q$-deformed
Conway-Coxeter friezes. We would have loved to have had the chance to talk all of this over with John.

\bigskip

%\begin{center}
%\includegraphics[scale=0.15]{photoJSV}
%\end{center}

%\newpage

%%%%%%%%%%%%%%%%%%%%
%%%%%%%%%%%%%%%%%%%%
\section{$q$-deformed Conway-Coxeter friezes}\label{CoCoSec}
%%%%%%%%%%%%%%%%%%%%
%%%%%%%%%%%%%%%%%%%%

Friezes, or frieze patterns, are perhaps one of the most ingenious and the least known of Coxeter's inventions~\cite{Cox}.  As he confessed, friezes
``caused [him] many restless nights,''~\cite{CoR} but Conway~\cite{CoCo} was able to put his mind somewhat to rest in answering his question about
the classification of friezes of positive integers.  As a result, these friezes are often called Conway-Coxeter friezes. Here, we introduce the
notion of $q$-deformations of Conway-Coxeter friezes.

%%%%%%%%%%%%%%%%%%%%
\subsection{$q$-analogues: why, what, where?}
%%%%%%%%%%%%%%%%%%%%
In mathematics and theoretical physics, ``$q$-deformation'' often means ``quantization'' and vice-versa. Historically, $q$ is the exponential of the
Planck constant $q=e^{\hbar},$ but, in quantization theory, both $q$~and~$\hbar$ are parameters.

A number of $q$-deformations of algebraic, geometric, and analytic structures have been introduced and thoroughly studied in mathematics. Among the
$q$-deformed structures, we encounter quantum groups, quantized Poisson structures, $q$-deformed special functions, just to mention the best known
theories. Quantized sequences of integers arise and play an important role in all of them.

Quantized quantities are functions in~$q$, usually polynomials, or power series.
A ``good'' $q$-analogue of a ``classical'' quantity must satisfy (at least) two requirements:
\begin{enumerate}
\item
when $q\to1$, we obtain the initial quantity;
\item
coefficients of polynomials in the quantized object  have a combinatorial meaning.
\end{enumerate}

\noindent Whenever a combinatorial object counts something, its $q$-analogue counts the same things but with more precision. This is why
combinatorics is playing an increasingly important role in mathematical physics.

%%%%%%%%%%%%%%%%%%%%
\subsection{Euler and Gauss}
%%%%%%%%%%%%%%%%%%%%

Quantum integers appeared in mathematics long before the development of quantum physics. For a positive integer~$n$, the polynomial
\begin{equation}
\label{Euler}
[n]_q:=
1+q+q^2+\cdots+q^{n-1}=\frac{1-q^n}{1-q}
\end{equation}
is commonly called the $q$-analogue of~$n$.
It satisfies two recurrences:
\begin{equation}
\label{RecEuler}
\left[n+1\right]_{q}=
q\left[n\right]_{q}+1,
\qquad\qquad
\left[n+1\right]_{q}=\left[n\right]_{q}+q^n,
\end{equation}
both of which can be used to define $[n]_q$ starting from the natural assumption $[0]_q=0$. The $q$-analogues~\eqref{Euler} were introduced by Euler,
who studied the infinite product $\prod_{n\geq1}(1-q^n)$, now called the Euler function. Euler connected this function to permutations simultaneously
founding combinatorics and the theory of modular forms.

The $q$-factorial is defined by~$[n]_q!:=[1]_q[2]_q\cdots[n]_q$,
and the Gaussian $q$-{\it binomial coefficients} by
$$
{n\choose m}_q:=
\frac{[n]_q!}{[m]_q!\,[n-m]_q!}.
$$
They are polynomials.
The role of $q$-binomials in combinatorics is immense.
They count points in Grassmannians over finite fields, Young diagrams,
binary words, etc.
Every coefficient of~${n\choose m}_q$ has a combinatorial meaning.

The first interesting example of a $q$-binomial is ${4\choose 2}_q=1+q+2q^2+q^3+q^4$ which is different from~$[6]_q$. Moreover,
$[3]_q!=1+2q+2q^2+q^3$ gives yet another version of ``quantum~$6$.'' These and many other examples explain our viewpoint: we cannot quantize~$6$ or
any other integer individually. What we quantize are not integers, but {\it sequences of integers}.

%%%%%%%%%%%%%%%%%%%%
\subsection{Friezes and Conway's classification}
%%%%%%%%%%%%%%%%%%%%

A frieze is an array of $n$ staggered infinite rows of positive integers, with the first and the last rows consisting of $1$'s, that satisfies the
local unimodular rule
$$
\xymatrix @!0 @R=0.45cm @C=0.45cm
{
\cdots&&1&& 1&&1&&1&&1&&\cdots&&1&&\cdots
 \\
&\cdots&&c_{i}\ar@{.}[rrrrrddd]&&c_{i+1}&&c_{i+2}&&\cdots&&c_{j-1}&&c_j\ar@{.}[lllllddd]&&\cdots
\\
 \\
 \\
&&&&&&&&c_{ij}\\
&&&&&&\cdots&&\cdots&&\cdots\\
\cdots&&1&& 1&&1&&1&&1&&\cdots&&1&&\cdots
}
\qquad\qquad
\xymatrix @!0 @R=0.45cm @C=0.45cm
{
\\
&\hbox{local rule}\\
&b\ar@{-}[rd]\ar@{-}[ld]&\\
a&&d\\
&c\ar@{-}[ru]\ar@{-}[lu]&\\
&ad-bc=1\\
}
$$
A frieze is determined by the second row $(c_i)_{i\in\Z}$. Every entry of the frieze is parametrized by two integers, $i$ and $j$. The local rule
then allows for the calculation of $c_{ij}$. For instance, $c_{ii}=c_i$, $c_{i,i+1}=c_ic_{i+1}-1$, etc.

\begin{ex}
\label{FEx}
The following $7$-periodic frieze with $6$ rows is the main example of~\cite{CoCo}
(and the favorite example of Coxeter who used it in all of his articles on the subject).
$$
\xymatrix @!0 @R=0.45cm @C=0.45cm
{
\cdots&&1&& 1&&1&&1&&1&&1&&1&&1&&1
 \\
&&&1&&4&&2&&1&&3&&2&&2&&1&&\cdots
 \\
\cdots&&1&&3&&7&&1&&2&&5&&3&&1
\\
&1&&2&&5&&3&&1&&3&&7&&1&&\cdots
 \\
 \cdots&&1&&3&&2&&2&&1&&4&&2&&1
 \\
&1&& 1&&1&&1&&1&&1&&1&&1&&1&&\cdots
}
\qquad\qquad
\xymatrix @!0 @R=0.35cm @C=0.35cm
{
&&&1\ar@{-}[rrd]\ar@{-}[lld]
\\
&3\ar@{-}[ldd]\ar@{-}[rrrr]&&&& 2\ar@{-}[rdd]&\\
\\
2\ar@{-}[rdd]&&&&&& 4\ar@{-}[ldd]\ar@{-}[llllldd]\ar@{-}[llllll]\ar@{-}[llllluu]\\
\\
&2\ar@{-}[rrrr]&&&&1
}
$$
\end{ex}

Coxeter proved~\cite{Cox} that every frieze is $(n+1)$-periodic. The Conway and Coxeter theorem provides a classification of friezes, stating that
friezes are in one-to-one correspondence with triangulations of a convex $(n+1)$-gon.

\begin{thm}[\cite{CoCo}]
\label{CoCoT} A sequence of positive integers $(c_0,\ldots,c_n)$ is a cycle of the second row of an $n$-row frieze if and only if $c_i$ is the number
of triangles at the $i$-th vertex of a triangulated $(n+1)$-gon.
\end{thm}

\noindent As Coxeter~\cite{CoR} acknowledged, this statement is actually due to John Conway.

Coxeter showed a connection of friezes to continued fractions; see Proposition~\ref{CoxPro} below. This connection is somewhat equivalent to the
existence of a natural embedding of friezes and triangulated $(n+1)$-gons into the Farey graph; see Figure~\ref{wtFg}.

Interest in Conway-Coxeter friezes has much increased recently because of the connection to various areas of number theory, algebra, geometry and
combinatorics. For a survey, see~\cite{Sop} and~\cite{Bau}.

%%%%%%%%%%%%%%%%%%%%
\subsection{Quantum friezes}
%%%%%%%%%%%%%%%%%%%%
The following notion is new, although it was implicitly in~\cite{SVRat}.

\begin{defn}
\label{qFrieze} A $q$-{\it deformed frieze} is an array of $n$ infinite rows of polynomials in one variable, with the first row of~$1$'s and the
second row of Euler's $q$-integers~\eqref{Euler},
$$
\xymatrix @!0 @R=0.5cm @C=0.65cm
{
\cdots&&1&& 1&&1&&1&&1&&\cdots&&1&&\cdots
 \\
&\cdots&&{\left[c_{i}\right]_q}\ar@{.}[rrrrrddd]&&{\left[c_{i+1}\right]_q}
&&{\left[c_{i+2}\right]_q}&&\cdots&&{\left[c_{j-1}\right]_q}&&\left[c_j\right]_q\ar@{.}[lllllddd]&&\cdots
\\
 \\
 \\
&&&&&&&&C_{ij}(q)\\
&&&&&&\cdots&&\cdots&&\cdots
}
$$
that satisfies the following ``$q$-unimodular rule''
\begin{equation}
\label{QL}
C_{i,j-1}(q)\,C_{i+1,j}(q)-C_{i+1,j-1}(q)\,C_{ij}(q)
=
q^{\displaystyle\sum_{k=i}^{j-1}(c_k-1)}.
\end{equation}
%We encircled the entries in the first row that contribute in the right-hand-side of~\eqref{QL}.
\end{defn}

Starting from the line $(\left[c_{i}\right]_q)$ and using~\eqref{QL},
one can calculate every next row of the $q$-frieze inductively.

\begin{ex}
\label{qFEx} The $q$-deformation of the frieze in Example~\ref{FEx} is
$$
\xymatrix @!0 @R=0.5cm @C=0.7cm
{
\cdots&&1&& 1&&1&&1&&1&&1&&1&&1&&1
 \\
&&&1&&[4]_q&&[2]_q&&1&&[3]_q&&[2]_q&&[2]_q&&1&&\cdots
 \\
\cdots&&1&&q[3]_q&&\{7\}_{q}&&1&&q[2]_q&&\{5\}_{q}&&[3]_{q}&&1
\\
&1&&q^2[2]_q&&q\{5\}_{q}&&[3]_q&&q^2&&q[3]_q&&\{7\}_{q}&&1&&\cdots
 \\
 \cdots&&q^3&&q^2[3]_q&&q[2]_q&&q^2[2]_q&&q^3&&q[4]_q&&[2]_q&&q^3
 \\
&q^3&&q^4&& q^2&&q^3&&q^3&&q^4&&q&&q^3&&q^4&&\cdots
}
$$
where $\{7\}_{q}=1+2q+2q^2+q^3+q^4$ and $\{5\}_{q}=1+2q+q^3+q^4$, and $[c]_q$ is as in~\eqref{Euler}. For instance, $\{7\}_{q}$ is calculated twice,
the first time as $\{7\}_{q}=[4]_q[2]_q-q^3$ and then as $\{7\}_{q}=\frac{\{5\}_{q}[3]_q-q^3}{[2]_q}$.
\end{ex}

At first glance, definition~\ref{qFrieze} does not look natural. In particular,~\eqref{QL} is no longer local, but it leads to a frieze with nice
properties. We will see in Section~\ref{ModSec}, that $C_{ij}(q)$ is a polynomial with positive integer coefficients.  In particular, the last row of
a $q$-frieze consists of powers of~$q$. Quantum friezes deserve thorough study, and we are currently at work on this project.

%%%%%%%%%%%%%%%%%%%%
%%%%%%%%%%%%%%%%%%%%
\section{Introducing quantum rationals \& irrationals}\label{IntSec}
%%%%%%%%%%%%%%%%%%%%
%%%%%%%%%%%%%%%%%%%%
Like the integers, rational numbers~$\frac{r}{s}$ cannot be $q$-deformed on their own, without including them in a sequence. Attempts to $q$-deform
the numerator and denominator separately lead to notions that lack nice properties. A very naive formula
$\left[\frac{r}{s}\right]_{q}=\frac{[r]_q}{[s]_q}$ and a more standard one $\left[x\right]_{q}=\frac{1-q^x}{1-q}$ (note that they coincide modulo a
rescaling of the parameter~$q$) are among them.
%The latter formula satisfies better recurrences, but sadly, is not a rational function when $x\in\Q$.

Here we give several equivalent definitions of $q$-rationals and explain the connection to $q$-friezes. For more definitions, a combinatorial
interpretation, and a connection to Jones polynomial and cluster algebras, see~\cite{SVRat}. Surprisingly, the $q$-deformation $[x]_q$ of an
irrational~$x\in\R$ is a Laurent series in~$q,$ owing to an unexpected phenomenon of stabilization of Taylor series of sequences of $q$-deformed
rationals converging to~$x$.

%%%%%%%%%%%%%%%%%%%%
\subsection{Deformed continued fractions}
%%%%%%%%%%%%%%%%%%%%

Every rational number~$\frac{r}{s}>0$ , where $r,s\in\Z_{>0}$ are coprime, has a standard finite continued fraction expansion
$\frac{r}{s}=[a_1,a_2,\ldots]$. Choosing an even number of coefficients (and removing the ambiguity $[a_1,\ldots,a_{n},1]=[a_1,\ldots,a_{n}+1]$), we
have the unique expansion $\frac{r}{s}=[a_1,\ldots,a_{2m}]$. Similarly, there is also a unique expansion with minus signs called the Hirzebruch-Jung
continued fraction:
$$
\frac{r}{s}
\quad=\quad
a_1 + \cfrac{1}{a_2
          + \cfrac{1}{\ddots +\cfrac{1}{a_{2m}} } }
           \quad =\quad
c_1 - \cfrac{1}{c_2
          - \cfrac{1}{\ddots - \cfrac{1}{c_k} } } ,
$$
where $a_i\geq1$ and $c_j\geq2$ (except for $a_1\geq0,c_1\geq1$).
The notation used by Hirzebruch is $\frac{r}{s}=\llbracket{}c_1,\ldots,c_k\rrbracket{}$;
the coefficients $a_i$ and $c_j$
are connected by the Hirzebruch
formula; see, e.g.,~\cite{MGO}.

\begin{defn}
\label{DefRat}
The $q$-deformed regular continued fraction is defined by
\begin{equation}
\label{qa}
[a_{1}, \ldots, a_{2m}]_{q}:=
[a_1]_{q} + \cfrac{q^{a_{1}}}{[a_2]_{q^{-1}}
          + \cfrac{q^{-a_{2}}}{[a_{3}]_{q}
          +\cfrac{q^{a_{3}}}{[a_{4}]_{q^{-1}}
          + \cfrac{q^{-a_{4}}}{
        \cfrac{\ddots}{[a_{2m-1}]_q+\cfrac{q^{a_{2m-1}}}{[a_{2m}]_{q^{-1}}}}}
          } }}
\end{equation}
where~$[a]_q$ is the Euler $q$-integer.
The $q$-deformed Hirzebruch-Jung continued fraction is
\begin{equation}
\label{qc}
\llbracket{}c_1,\ldots,c_k\rrbracket_{q}:=
[c_1]_{q} - \cfrac{q^{c_{1}-1}}{[c_2]_{q}
          - \cfrac{q^{c_{2}-1}}{\ddots \cfrac{\ddots}{[c_{k-1}]_{q}- \cfrac{q^{c_{k-1}-1}}{[c_k]_{q}} } }}
\end{equation}
For $\frac{r}{s}=[a_1,\ldots,a_{2m}]=\llbracket{}c_1,\ldots,c_k\rrbracket{}$, the rational functions~\eqref{qa} and~\eqref{qc} coincide. This
rational function is called the $q$-{\it rational} and denoted~$\left[\frac{r}{s}\right]_{q}=\frac{\Rc(q)}{\Sc(q)}.$
\end{defn}

Both polynomials,~$\Rc$ and~$\Sc$, depend on~$r$ and~$s$.

\begin{ex}
\label{Q5}
For instance, we obtain
$\left[\frac{5}{2}\right]_{q}=
\frac{1+2q+q^{2}+q^{3}}{1+q}$
and
$\left[\frac{5}{3}\right]_{q}=
\frac{1+q+2q^{2}+q^{3}}{1+q+q^{2}}.$
Observe that ``quantum~$5$'' in the numerator depends on the denominator.
\end{ex}

%%%%%%%%%%%%%%%%%%%%
\subsection{The weighted Farey graph and Stern-Brocot tree}\label{FaSec}
%%%%%%%%%%%%%%%%%%%%
Let us give a recursive definition.

The set of rational numbers $\Q$, completed by $\infty:=\frac{1}{0}$, are vertices of a graph called the {\it Farey graph}. Two
rationals~$\frac{r}{s}$ and~$\frac{r'}{s'}$ are connected by an edge if and only if $rs'-r's=\pm1$. Edges of the Farey graph are often represented as
geodesics of the hyperbolic plane which is triangulated.

\begin{figure}
\begin{center}
\includegraphics[width=10cm,height=5.5cm]{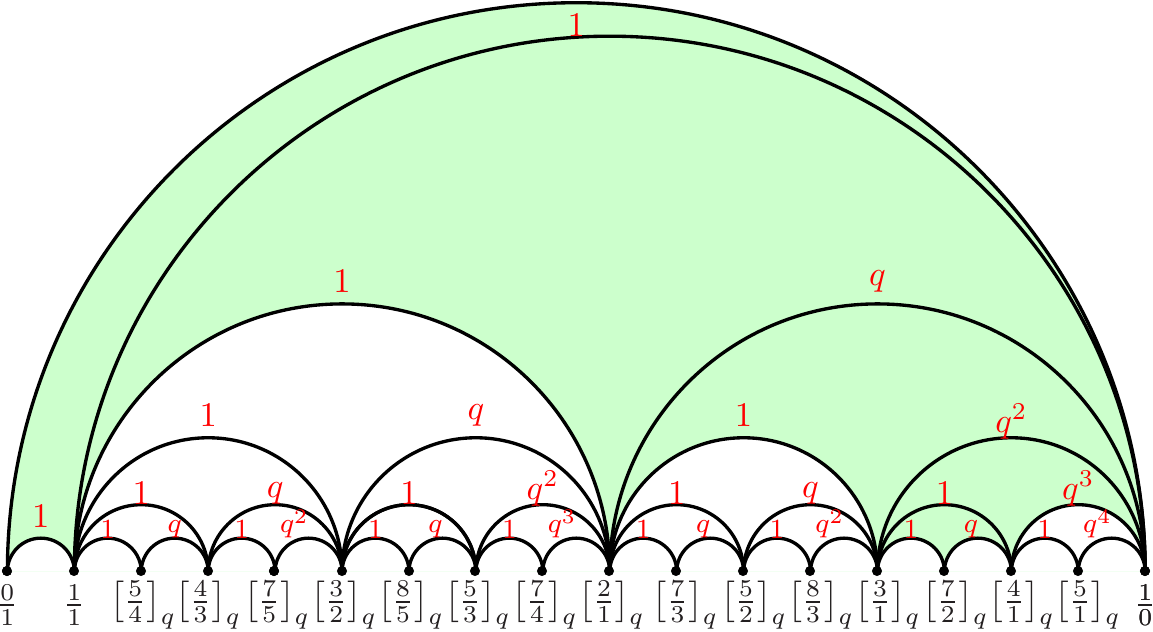}
\caption{The weighted Farey graph between~$\frac01$ and~$\frac10$; the colored area corresponds to
the frieze of Example~\ref{FEx}.}
\label{wtFg}
\end{center}
\end{figure}

\begin{defn}
\label{RecDef} The {\it weighted Farey graph}, see Figure~\ref{wtFg}, is just the classical Farey graph, in which the vertices are labeled by
rational functions in~$q$ and the edges are weighted by powers of~$q$. The weights and the labels are defined recursively via  the following local
rule:
\begin{center}
\includegraphics[width=3cm,height=2.3cm]{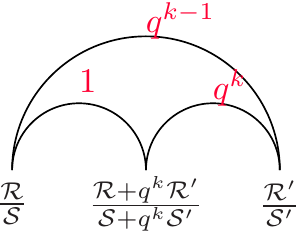}
\end{center}
from the initial triangle $(\frac{0}{1},\frac{1}{1},\frac{1}{0})$
that remains undeformed.
\end{defn}

\begin{thm}[\cite{SVRat}]
The vertices of the weighted Farey graph are, indeed, labeled by $q$-rationals.
\end{thm}

\begin{rem}
The Farey graph contains~$\Z$ as a subgraph (that forms a sequence of triangles, $(\frac{n}{1},\frac{n+1}{1},\frac{1}{0})$). Our $q$-deformation
restricted to~$\Z$ leads to Euler's formula~\eqref{Euler}.
\end{rem}

Alternatively, one can use the Stern-Brocot tree instead of the Farey graph.
$$
\xymatrix @!0 @R=0.43cm @C=0.35cm
{
&&&&&&&&&&&&&&&&\ar@{-}[dd]\\
&&&&&&&&&&&&&&\frac{1}{0}&&&&\frac{0}{1}\\
&&&&&&&&&&&&&&&&\bullet\ar@{-}[lllllllldd]_{\textcolor{red}{1}}\ar@{-}[rrrrrrrrdd]^{\textcolor{red}{1}}&&&&&&&&\\
&&&&&&&&&&&&&&&&\frac{1}{1}\\
&&&&\frac{1}{0}&&&&\bullet\ar@{-}[lllldd]_{\textcolor{red}{q}}\ar@{-}[rrrrdd]^{\textcolor{red}{1}}
&&&&&&&&&&&&&&&&\bullet\ar@{-}[lllldd]_{\textcolor{red}{q}}\ar@{-}[rrrrdd]^{\textcolor{red}{1}}&&&&\frac{0}{1}\\
&&&&&&&&\left[\frac{2}{1}\right]_q
&&&&&&&&&&&&&&&&\left[\frac{1}{2}\right]_q\\
&&&&\bullet\ar@{-}[lldd]_{\textcolor{red}{q^2}}\ar@{-}[rrdd]^{\textcolor{red}{1}}
&&&&&&&&\bullet\ar@{-}[lldd]_{\textcolor{red}{q}}\ar@{-}[rrdd]^{\textcolor{red}{1}}
&&&&&&&&\bullet\ar@{-}[lldd]_{\textcolor{red}{q^2}}\ar@{-}[rrdd]^{\textcolor{red}{1}}
&&&&&&&&\bullet\ar@{-}[lldd]_{\textcolor{red}{q}}\ar@{-}[rrdd]^{\textcolor{red}{1}}\\
&&&&&&&&&&&&&&&&&&&&&&&&&&&&\\
&&\bullet\ar@{-}[ldd]_{\textcolor{red}{q^3}}\ar@{-}[rdd]^{\textcolor{red}{1}}
&&\;\left[\frac{3}{1}\right]_q&&\bullet\ar@{-}[ldd]_{\textcolor{red}{q}}\ar@{-}[rdd]^{\textcolor{red}{1}}
&&&&\bullet\ar@{-}[ldd]_{\textcolor{red}{q^2}}\ar@{-}[rdd]^{\textcolor{red}{1}}
&&\;\left[\frac{3}{2}\right]_q&&\bullet\ar@{-}[ldd]_{\textcolor{red}{q}}\ar@{-}[rdd]^{\textcolor{red}{1}}
&&&&\bullet\ar@{-}[ldd]_{\textcolor{red}{q^3}}\ar@{-}[rdd]^{\textcolor{red}{1}}
&&\;\left[\frac{2}{3}\right]_q&&\bullet\ar@{-}[ldd]_{\textcolor{red}{q}}\ar@{-}[rdd]^{\textcolor{red}{1}}
&&&&\bullet\ar@{-}[ldd]_{\textcolor{red}{q^2}}\ar@{-}[rdd]^{\textcolor{red}{1}}
&&\;\left[\frac{1}{3}\right]_q&&\bullet\ar@{-}[ldd]_{\textcolor{red}{q}}\ar@{-}[rdd]^{\textcolor{red}{1}}\\
&&&&&&&&&&&&&&&&&&&&&&
&&&&&&&&\\
&&\;\left[\frac{4}{1}\right]_q
&&&&\;\left[\frac{5}{2}\right]_q
&&&&\;\left[\frac{5}{3}\right]_q
&&&&\;\left[\frac{4}{3}\right]_q
&&&&\;\left[\frac{3}{4}\right]_q
&&&&\;\left[\frac{3}{5}\right]_q
&&&&\;\left[\frac{2}{5}\right]_q&&&&\;\left[\frac{1}{4}\right]_q&&\\
&&&&\ldots&&&&&&&&&&&&\ldots&&&&&&&&&&&&\ldots
}
\qquad
\xymatrix @!0 @R=0.45cm @C=0.4cm
{
\\
\\
\\
&&&\mathrm{(local}\; \mathrm{rule)}\\
&&&\ar@{-}[ddd]^{\textcolor{red}{q^{k-1}}}\\
\\
&\frac{\Rc}{\Sc}&&&&&\!\!\!\!\!\frac{\Rc'}{\Sc'}\\
&&&\bullet\ar@{-}[llldd]_{\textcolor{red}{q^k}}\ar@{-}[rrrdd]^{\textcolor{red}{1}}\\
\\
&&&\frac{q^k\Rc+\Rc'}{q^k\Sc+\Sc'}&&&&&&&&
}
$$
The weight of every edge of the weighted Stern-Brocot tree is determined by the local rule along the tree: the weight of the right branch is
multiplied by~$q$, the weight of the left branch is~$1$. This allows us to calculate $q$-rationals inductively. For instance, $
\left[\frac{2}{1}\right]_q= \frac{q+1}{1}=[2]_q, \; \left[\frac{1}{2}\right]_q= \frac{q}{q+1}. $ The left branch of the tree consists of the
classical $q$-integers~\eqref{Euler}.

%%%%%%%%%%%%%%%%%%%%
\subsection{$q$-deformed $\SL(2,\Z)$}\label{ModSec}
%%%%%%%%%%%%%%%%%%%%

The {\it modular group} $\SL(2,\Z)$ is useful as a tool in working with continued fractions. In our case, we need to $q$-deform $\SL(2,\Z)$, so we
consider the following three matrices:
\begin{equation}
\label{SL2q}
R_{q}:=
\begin{pmatrix}
q&1\\[2pt]
0&1
\end{pmatrix},
\quad
L_{q}:=
\begin{pmatrix}
q&0\\[2pt]
q&1
\end{pmatrix},
\quad
S_{q}:=
\begin{pmatrix}
0&-q^{-1}\\[4pt]
1&0
\end{pmatrix}.
\end{equation}
For $q=1$, we obtain the standard matrices, any two of which can be chosen as generators of~$\SL(2,\Z)$.

\begin{prop}[\cite{SVRat}]
\label{SLProP} Given a rational, $\frac{r}{s}=[a_1,\ldots,a_{2m}]=\llbracket{}c_1,\ldots,c_k\rrbracket$, the polynomials $\Rc(q)$ and~$\Sc(q)$
of~$\left[\frac{r}{s}\right]_q=\frac{\Rc(q)}{\Sc(q)}$ can be calculated as the entries of the first column of the following matrix products
$$
R_{q}^{a_{1}}L_{q}^{a_{2}}\cdots R_{q}^{a_{2m-1}}L_{q}^{a_{2m}}=
\begin{pmatrix}
q\Rc&\Rc'\\[2pt]
q\Sc&\Sc'
\end{pmatrix},
\qquad\qquad
R_{q}^{c_{1}}S_{q}R_{q}^{c_{2}}S_{q} \cdots R_{q}^{c_{k}}S_{q}=
\begin{pmatrix}
\Rc&-q^{c_{k}-1}\Rc''\\[6pt]
\Sc&-q^{c_{k}-1}\Sc''
\end{pmatrix},
$$
where $\Rc',\Sc'$ and $\Rc'',\Sc''$ are lower degree polynomials corresponding to
the previous convergents of the continued fractions.
\end{prop}

This way of understanding $q$-rationals clarifies Definition~\ref{qFrieze}. Indeed,~\eqref{QL} is the determinant of the second matrix. Note that the
matrices~\eqref{SL2q} are proportional to the matrices arising in quantum Teichm\"uller theory; see~\cite{CF}. This interesting connection has yet to
be investigated.

%%%%%%%%%%%%%%%%%%%%
\subsection{Connection to $q$-friezes}
%%%%%%%%%%%%%%%%%%%%

It turns out that $q$-rationals appear everywhere in $q$-deformed Conway-Coxeter friezes,
as quotients of the neighboring entries, so that
we can consider $q$-friezes as yet another way to describe $q$-deformations of rationals.

Let us first recall the connection between the Conway-Coxeter friezes and
continued fractions.

\begin{prop}[\cite{Cox}]
\label{CoxPro}
If $(c_{ij})$ are the entries of a frieze, then
$\frac{c_{ij}}{c_{i+1,j}}=\llbracket{}c_i,\ldots,c_j\rrbracket$.
\end{prop}

This statement has a straightforward $q$-analogue.

\begin{prop}
\label{FQThm}
If $(c_{ij})$ are the entries of a frieze and $(C_{ij}(q))$ the entries of the $q$-deformed frieze, then
%For every $i,j$, such that $|i-j|\leq n-2$, one has
$
\frac{C_{ij}(q)}{C_{i+1,j}(q)}=\left[\frac{c_{ij}}{c_{i+1,j}}\right]_q=\llbracket{}c_i,\ldots,c_j\rrbracket_q.
$
\end{prop}

The $q$-frieze of Example~\ref{qFEx} contains many examples of $q$-rationals:
$\frac{\{5\}}{[2]}=\left[\frac{5}{2}\right]_q$, $\frac{\{7\}}{[3]}=\left[\frac{7}{3}\right]_q$, etc.

%%%%%%%%%%%%%%%%%%%%
\subsection{$q$-deformed irrational numbers}\label{StabSec}
%%%%%%%%%%%%%%%%%%%%

Let $x\geq0$ be an irrational number, and~let $(x_n)_{n\geq1}$ be a sequence of rationals converging to~$x$. Our definition is the following.

Take the sequence of rational functions $\left[x_1\right]_q,\left[x_2\right]_q,\ldots$ Consider their Taylor expansions at~$q=0$, for which we will
use the notation $\left[x_n\right]_q=\sum_{k\geq0}\varkappa_{n,k}\,q^k.$

\begin{thm}[\cite{SVRe}]
\label{ConvThm}
(i) For every~$k\geq0$, the coefficients of the Taylor series of~$\left[x_n\right]_q$ stabilize as~$n$ grows.

(ii) The limit coefficients~$\varkappa_k=\lim_{n\to\infty}\varkappa_{n,k}$
do not depend on the sequence $(x_n)_{n\geq1}$,
but only on~$x$.
\end{thm}

This stabilization phenomenon allows us to define the $q$-deformation of~$x\geq0$ as a power series in~$q$:
$$
\left[x\right]_q=\varkappa_0+\varkappa_1q+\varkappa_2q^2+\varkappa_3q^3+\cdots
$$
In the case of~$x<0$, the definition of $q$-deformation is based on the recurrence $\left[x-1\right]_q:=q^{-1}\left[x\right]_q-q^{-1}$; cf.
Section~\ref{SLAct}. It turns out that the resulting series in~$q$ is a {\it Laurent series} (with integer coefficients):
$$
\left[x\right]_q=-q^{-N}+\varkappa_{1-N}\,q^{1-N}+\varkappa_{2-N}\,q^{2-N}+\cdots ,
$$
where $N\in\Z_{>0}$ such that $-N\leq{}x<1-N$.

%%%%%%%%%%%%%%%%%%%%
%%%%%%%%%%%%%%%%%%%%
\section{Examples: $q$-Fibonacci and $q$-Pell numbers,
$q$-golden ratio and~$\left[\sqrt{2}\right]_q$}\label{ExSec}
%%%%%%%%%%%%%%%%%%%%
%%%%%%%%%%%%%%%%%%%%

Let us consider two remarkable sequences of rationals,
$$
\frac{F_{n+1}}{F_{n}}=
\underbrace{\left[1, 1, \ldots,1\right]}_n
\qquad\hbox{and}\qquad
\frac{P_{n+1}}{P_{n}}=
\underbrace{\left[2, 2, \ldots,2\right]}_n,
$$
where~$F_n$ is the~$n^{\hbox{th}}$ Fibonacci number,
and~$P_n$ is the~$n^{\hbox{th}}$ Pell number.
Quantizing them, we obtain sequences of polynomials
with a very particular ``swivel'' property:
the polynomials corresponding to~$F_n$ (resp.~$P_n$) in the numerator and denominator
are mirrors of each other.
The stabilized Taylor series of $\left[\frac{F_{n+1}}{F_{n}}\right]_q$ and $\left[\frac{P_{n+1}}{P_{n}}\right]_q$
give rise of $q$-analogues of $\frac{1+\sqrt{5}}{2}$ and~$1+\sqrt{2}$, respectively.

%%%%%%%%%%%%%%%%%%%%
\subsection{$q$-Fibonacci numbers}\label{GoR}
%%%%%%%%%%%%%%%%%%%%

Let~$\F_n(q)$ be a sequence of polynomials defined by the recurrence
$$
\F_{n+2}= [3]_q\,\F_n-q^2\F_{n-2},
$$
where $[3]_q=1+q+q^2$ is  Euler's quantum~$3$, and where the initial conditions are $(\F_0(q)=0, \F_2(q)=1)$ and $(\F_1(q)=1, \F_3(q)=1+q). $ The
sequence of polynomials $\F_n(q)$ is a $q$-deformation of the Fibonacci sequence, i.e., $\F_n(1)=F_n$. Consider also the {\it mirror polynomials}
$\tilde\F_n(q):=q^{n-2}\F_n(\frac{1}{q})$ ($n\geq2$).
The triangles of their coefficients %of~$\F_n(q)$ and~$\tilde\F_n(q)$ (for~$n\geq2$) are
$$
\begin{array}{rcccccccc}
1\\
1&1\\
1&1&1\\
1&2&1&1\\
1&2&2&2&1\\
1&3&3&3&2&1\\
1&3&4&5&4&3&1\\
\cdots
\end{array}
\qquad\qquad
\begin{array}{rcccccccc}
&&&&&&1\\
&&&&&1&1\\
&&&&1&1&1\\
&&&1&1&2&1\\
&&1&2&2&2&1\\
&1&2&3&3&3&1\\
1&3&4&5&4&3&1\\
&&&&&\cdots
\end{array}
$$
are the well-studied sequences A079487 and  A123245 of OEIS~\cite{OEIS}.

\begin{prop}
\label{FiboT}
One has
$
\left[\frac{F_{n+1}}{F_{n}}\right]_q=
\frac{\tilde\F_{n+1}(q)}{\F_n(q)}.
$
\end{prop}

\begin{ex}
\label{FibpPEx} The case $\left[\frac{5}{3}\right]_q$ has already been considered in Example~\ref{Q5}. We then have
$$
\begin{array}{rcl}
\left[\frac{8}{5}\right]_q&=&\displaystyle\frac{1+2q+2q^2+2q^3+q^4}{1+2q+q^2+q^3},
\\[12pt]
\left[\frac{13}{8}\right]_q&=&\displaystyle\frac{1+2q+3q^2+3q^3+3q^4+q^5}{1+2q+2q^2+2q^3+q^4},
\\[12pt]
\left[\frac{21}{13}\right]_q&=&\displaystyle\frac{1+3q+4q^2+5q^3+4q^4+3q^5+q^6}{1+3q+3q^2+3q^3+2q^4+q^5},
\\
\ldots&\ldots&\ldots
\end{array}
$$
The coefficients of these rational functions grow at every fixed power of~$q$,
and there is no stabilization of rational functions.
\end{ex}

%%%%%%%%%%%%%%%%%%%%
\subsection{$q$-Pell numbers}\label{PeP}
%%%%%%%%%%%%%%%%%%%%

The sequence of polynomials~$\Pc_{n}(q)$ satisfying the recursion
$$
\Pc_{n+2} = {4\choose 2}_q\Pc_{n} - q^4\Pc_{n-2},
$$
where ${4\choose 2}_q$ is the Gaussian $q$-binomial and with the initial conditions $(\Pc_{0}(q)=0,\Pc_{2}(q)=1+q)$ and $(\Pc_{1}(q)=1,
\Pc_{3}(q)=1+q+2q^2+q^3)$, are $q$-analogues of the classical {\it Pell numbers}. The mirror polynomials are defined by
$\tilde\Pc_{n}:=q^{2n-3}\Pc_{n}(\frac{1}{q})$ ($n\geq2$). The coefficients of~$\Pc_{n}(q)$
$$
\begin{array}{rccccccccccc}
1\\
1&1\\
1&2&1&1\\
1&2&3&3&2&1\\
1&3&5&6&6&5&2&1\\
1&3&7&11&13&13&11&7&3&1\\
\cdots
\end{array}
$$
form a triangular sequence that was recently added to the OEIS (sequence A323670).

\begin{prop}
\label{PellT}
One has
$
\left[\frac{P_{n+1}}{P_{n}}\right]_q=
\frac{\tilde\Pc_{n+1}(q)}{\Pc_n(q)}.
$
\end{prop}

%%%%%%%%%%%%%%%%%%%%
\subsection{$q$-deformed golden ratio}
%%%%%%%%%%%%%%%%%%%%

To illustrate the stabilization phenomenon, we take the Taylor series of
the Fibonacci quotients; see Example~\ref{FibpPEx}.
For instance,
$$
\begin{array}{rcl}
\left[\frac{8}{5}\right]_q&=&
1 + q^2 - q^3 + 2 q^4 - 4 q^5 + 7 q^6 - 12 q^7 + 21 q^8 - 37q^9+ 65q^{10} - 114q^{11}+ 200q^{12}\cdots\\[6pt]
\left[\frac{21}{13}\right]_q&=&
1 + q^2 - q^3 + 2 q^4 - 4 q^5 + 8 q^6 - 17 q^7 + 36 q^8 - 75q^9+ 156q^{10} - 325q^{11} + 677q^{12}\cdots\\[6pt]
\left[\frac{55}{34}\right]_q&=&
1 + q^2 - q^3 + 2 q^4 - 4 q^5 + 8 q^6 - 17 q^7 + 37 q^8 - 82q^9+184q^{10} - 414q^{11}+ 932q^{12} \cdots\\[6pt]
\left[\frac{144}{89}\right]_q&=&
1 + q^2 - q^3 + 2 q^4 - 4 q^5 + 8 q^6 - 17 q^7 + 37 q^8 - 82 q^9 + 185 q^{10} - 423 q^{11} + 978 q^{12} \cdots
\end{array}
$$
More and more coefficients repeat as the series proceed with the series eventually stabilizing to
$$
\begin{array}{rcl}
\left[\varphi\right]_q&=&
1 + q^2 - q^3 + 2 q^4 - 4 q^5 + 8 q^6 - 17 q^7 + 37 q^8 - 82 q^9 + 185 q^{10} - 423 q^{11} + 978 q^{12}-2283q^{13}\\[4pt]
&&+ 5373q^{14}-12735q^{15}+30372q^{16}-72832q^{17}+175502q^{18}-424748q^{19}+1032004q^{20} \cdots
\end{array}
$$
This power series is our quantized golden ratio $\left[\frac{1+\sqrt{5}}{2}\right]_q:=\left[\varphi\right]_q$.

The series~$\left[\varphi\right]_q$ satisfies the equation
$
q\left[\varphi\right]^2_q-
\left(q^2+q-1 \right)\left[\varphi\right]_q -1 =0,
$
which is a $q$-analogue of $x^2-x+1=0$.
Therefore, the generating function of the series~$\left[\varphi\right]_q$ is
\begin{equation}
\label{Gold}
\left[\varphi\right]_q=
\frac{q^2+q-1+\sqrt{(q^2 + 3q + 1)(q^2 - q + 1)}}{2q}.
\end{equation}
Let us mention that the coefficients of the series $\left[\varphi\right]_q$
remarkably coincide with Sequence A004148 of~\cite{OEIS}
(with alternating signs), called the Generalized Catalan numbers.

%%%%%%%%%%%%%%%%%%%%
\subsection{$q$-deformed $\sqrt{2}$}\label{Qold2Sec}
%%%%%%%%%%%%%%%%%%%%

Quotients of the Pell polynomials stabilize to the series~$\left[1+\sqrt{2}\right]_q$,
from which we deduce
$$
\begin{array}{rcl}
\left[\sqrt{2}\right]_q&=&
1+ q^3 - 2q^5  + q^6+ 4q^7- 5q^8- 7q^9+18q^{10} + 7q^{11}- 55q^{12}+ 18q^{13}\\[4pt]
&&+ 146q^{14} - 155q^{15}- 322q^{16}+ 692q^{17}+ 476q^{18}
 - 2446q^{19}+ 307q^{20} \cdots .
 \end{array}
$$
This series is a solution of $ q^2\left[\sqrt{2}\right]^2_q-\left(q^3-1 \right)\left[\sqrt{2}\right]_q=q^2 + 1, $ which is our version of the
$q$-analogue of~$x^2=2$. The generating function is then equal to
\begin{equation}
\label{Silver}
\left[\sqrt{2}\right]_q=
\frac{q^3-1+\sqrt{(q^4+q^3+4q^2+q+1)(q^2-q+1)}}{2q^2}.
\end{equation}
Note that the coefficients of $\left[\sqrt{2}\right]_q$ grow much more slowly than those of~$\left[\varphi\right]_q$ and fail to match any known
sequence.

%%%%%%%%%%%%%%%%%%%%
\subsection{Quadratic $q$-irrationals}
%%%%%%%%%%%%%%%%%%%%
Several observations can be made by analyzing~\eqref{Gold} and~\eqref{Silver}. The polynomial under the radical of these $q$-numbers is a {\it
palindrome}. This remarkable property remains true for arbitrary {\it quadratic irrationals}, i.e., numbers of the form $x=\frac{a+\sqrt{b}}{c}$,
where $a,b>0,c$ are integers.

\begin{thm}[\cite{LMG}]
\label{QIThm}
For every quadratic irrational,
$
\left[\frac{a+\sqrt{b}}{c}\right]_q=
\frac{A(q)+\sqrt{B(q)}}{C(q)},
$
where $A$, $B$ and $C$ are polynomials in~$q$.
Furthermore,~$B$ is a monic polynomial whose coefficients form a palindrome.
\end{thm}

A quadratic irrational can also be characterized as a fixed point of an element of~$\PSL(2,\Z)$. A $q$-analogue of this property will be provided in
Theorem~\ref{recprop} in the next section.

%%%%%%%%%%%%%%%%%%%%
\subsection{Radius of convergence}
%%%%%%%%%%%%%%%%%%%%
The modulus of the smallest (i.e., closest to~$0$) root of the polynomials under the radical
in~\eqref{Gold} and~\eqref{Silver} is equal to
$$
R^-_\varphi=\frac{3-\sqrt{5}}{2}
\qquad\hbox{and}\qquad
R^-_{\sqrt{2}}=
\frac{1 + \sqrt{2} - \sqrt{2\sqrt{2} - 1}}{2},
$$
respectively. This are the {\it radii of convergence} of the Taylor series~\eqref{Gold} and~\eqref{Silver}, resp. The modulus of the largest roots
are $R^+_\varphi=1/R^-_\varphi$ and $R^+_{\sqrt{2}}=1/R^-_{\sqrt{2}}$, viz.,
$$
R^+_\varphi=\frac{3+\sqrt{5}}{2}
\qquad\hbox{and}\qquad
R^+_{\sqrt{2}}=
\frac{1 + \sqrt{2} + \sqrt{2\sqrt{2} - 1}}{2}.
$$

Aesthetic aspects of these formulas motivated us to analyze several more examples, one of which is another remarkable number,
$\frac{9+\sqrt{221}}{14}$, sometimes called the ``bronze ratio.'' This is the third, after $\varphi$ and~$\sqrt{2}$, badly approximated number in
Markov theory. The modulus of the minimal and maximal roots of the polynomial under the radical in $\left[\frac{9+\sqrt{221}}{14}\right]_q$ are
$$
R^-=
\frac{1 + \sqrt{13} - \sqrt{2\left(\sqrt{13} - 1\right)}}{4}
\qquad\hbox{and}\qquad
R^+=
\frac{1 + \sqrt{13} + \sqrt{2\left(\sqrt{13} - 1\right)}}{4}.
$$
Note that $221=13\cdot{}17$. We do not know if the striking resemblance to the case of~$\sqrt{2}$ is a coincidence. Work on the analytic properties
of $q$-numbers is in progress~\cite{Ves}.

%%%%%%%%%%%%%%%%%%%%
%%%%%%%%%%%%%%%%%%%%
\section{Some properties of $q$-rationals}\label{PropSec}
%%%%%%%%%%%%%%%%%%%%
%%%%%%%%%%%%%%%%%%%%

The notion of $q$-rationals arose from an attempt to understand the connection between several different theories, such as continued fractions, Jones
polynomials of (rational) knots, quantum Teichm\"uller theory, and cluster algebras. These connections were discovered by many authors;
see~\cite{CF,R2}. Our definition is a specialization of such notions as $F$-polynomials, quantum geodesic length, snake graphs. We present some of
those concrete properties of $q$-rationals that we consider most important; many reflect this deep connection. For more details,
see~\cite{SVRat,SVRe}.

%%%%%%%%%%%%%%%%%%%%
\subsection{$\PSL(2,\Z)$-action}\label{SLAct}
%%%%%%%%%%%%%%%%%%%%
Our first important property is the following.

\begin{thm}[\cite{LMG}]
\label{recprop}
The procedure of $q$-deformation commutes with the $\PSL(2,\Z)$-action.
\end{thm}
Indeed,
the first recurrence~\eqref{RecEuler} remains true for $q$-rationals.
For every~$x\in\Q$, we have
\begin{equation}
\label{TransEq}
\left[x+1\right]_q=q\left[x\right]_q+1.
\end{equation}
Recurrence~\eqref{TransEq} readily follows from~\eqref{qa} and is very useful as it allows us to define $q$-deformations of~$x<0$. Furthermore, we
have
\begin{equation}
\label{InvEq}
\left[-\frac{1}{x}\right]_q=-\frac{1}{q\left[x\right]_{q}},
\qquad\qquad
\left[-x\right]_q=-q^{-1}\left[x\right]_{q^{-1}}.
\end{equation}
Together,~\eqref{TransEq} and~\eqref{InvEq} define an action of~$\PSL(2,\Z)$ on $q$-rationals generated by the matrices~$R_q$ and~$qS_q$
in~\eqref{SL2q}. This was (implicitly) checked in~\cite{SVRat} (see Lemma~4.6) and will be further developed in~\cite{LMG}.

\begin{rem}
Recurrence~\eqref{TransEq} appears, e.g., for $q$-integers and plays a crucial role in quantum algebra. For instance, this recurrence is necessary
for the quantum binomial formula. Identities~\eqref{InvEq} look more intriguing and need to be better understood.
\end{rem}

%%%%%%%%%%%%%%%%%%%%
\subsection{Total positivity}
%%%%%%%%%%%%%%%%%%%%
The polynomials in the numerator and denominator of a positive $q$-rational
have positive integer coefficients, with~$1$ at lowest and higher orders.
Moreover,
the set of all $q$-rationals has a much stronger property of {\it total positivity}.
Consider two $q$-rationals, $\left[\frac{r}{s}\right]_{q}=\frac{\Rc(q)}{\Sc(q)}$
and $\left[\frac{r'}{s'}\right]_{q}=\frac{\Rc'(q)}{\Sc'(q)}$.

\begin{thm}[\cite{SVRat}]
\label{ToPoT}
If $\frac{r}{s}>\frac{r'}{s'}>0$, then the polynomial
\begin{equation}
\label{XPol}
\cX_{\frac{r}{s},\frac{r'}{s'}}(q):=\Rc(q)\Sc'(q)-\Sc(q)\Rc'(q)
\end{equation}
 has positive integer coefficients.
\end{thm}

The main ingredient of the proof of this theorem is the fact that~\eqref{XPol}
is a monomial, i.e., proportional to a power of~$q$, if and only if
$\frac{r}{s}$ and~$\frac{r'}{s'}$ are connected in the Farey graph.

Theorem~\ref{ToPoT} means that the ``quantization preserves the order,'' in the sense that it is a homeomorphism of~$\Q$ into an ordered subset of a
partially ordered set of rational functions. The notion of total positivity has a long history in mathematics and manifests itself in every area of
it.

%%%%%%%%%%%%%%%%%%%%
\subsection{Relation to the Jones polynomial}
%%%%%%%%%%%%%%%%%%%%

One of the ingenious inventions of Conway~\cite{Con} was to encode a certain large class of knots, called ``rational'' or ``two-bridge'' knots, by
continued fractions. Every notion of knot theory, such as knot invariants, can then be directly associated to continued fractions.

The Jones polynomial is a powerful invariant in knot theory. It turns out that one can express the Jones polynomial $J_{\frac{r}{s}}(q)$ of the
two-bridge knot associated to a rational $\frac{r}{s}$ as a combination of the polynomials $\Rc$ and~$\Sc$ of the $q$-rational
$\left[\frac{r}{s}\right]_q$.

\begin{thm}[\cite{SVRat}]
\label{JRS}
The Jones polynomial of a two-bridge knot is
$
J_{\frac{r}{s}}(q)=q\Rc(q)+(1-q)\Sc(q).
$
\end{thm}
\noindent The proof is based on the connection to cluster algebras; see~\cite{R2}.

%%%%%%%%%%%%%%%%%%%%
\subsection{Unimodality conjecture}
%%%%%%%%%%%%%%%%%%%%
Long computer experimentation and evidence in the simplest cases convinced us of yet another property of $q$-rationals.

\begin{conj}[\cite{SVRat}]
For every $\left[\frac{r}{s}\right]_{q}=\frac{\Rc(q)}{\Sc(q)}$,
the coefficients of the polynomials~$\Rc(q)$ and~$\Sc(q)$ form unimodal sequences.
\end{conj}

This means that the coefficients of the polynomials increase from~$1$ to the maximal value
(that can be taken by one or more consecutive coefficients) and then decrease to~$1$.
Unimodal sequences appear in mathematics, and this property is interesting because
it often hides some combinatorial or geometric structure.

A proof of several particular cases of the conjecture, as well a connection to old problems of combinatorics, was obtained in~\cite{CSS}. However,
the general problem is open.

%%%%%%%%%%%%%%%%%%%%
%%%%%%%%%%%%%%%%%%%%
\section{Discussion Actually}\label{Conclusion}
%%%%%%%%%%%%%%%%%%%%
%%%%%%%%%%%%%%%%%%%%
To our great regret, we cannot discuss all of this with John Conway (at least not yet!), but we can easily imagine such a discussion.

Perhaps, John would ask about {\it number walls}, as he did when we talked about friezes. This notion, among hundreds of other surprising notions,
can be found on leafing through John's co-authored (with Richard Guy) book~\cite{CG}. A number wall is a~$\Z^2$-lattice filled with integers
satisfying a local rule different from the frieze rule:
$$
\xymatrix @!0 @R=0.55cm @C=0.55cm
{
\cdots&\bullet&N&\bullet&\cdots\\
\cdots&W&X&E&\cdots\\
\cdots&\bullet&S&\bullet&\cdots
}
\qquad\qquad
\xymatrix @!0 @R=0.45cm @C=0.45cm
{
\\
\hbox{local rule}\\
X^2=NS+EW.\\
}
$$
As in the case of friezes, the array is bounded by lines (or zig-zags) of~$1$'s. ``Can you $q$-deform this?,'' John would likely ask. Our answer
would be: ``We don't know!'' The question, however, is interesting since number walls produce many interesting (known and unknown) sequences of
integers.

He might also have wondered about {\it surreal numbers}, which he and Guy also discussed in~\cite{CG}. Surreal numbers can be constructed using the
same binary tree as in Section~\ref{FaSec}. Vertices of the tree are labeled by dyadic rationals (rationals whose denominator is a power of~$2$):
$$
\xymatrix @!0 @R=0.43cm @C=0.35cm
{
&&&&&&&&&&&&&&&&0\ar@{-}[lllllllldd]\ar@{-}[rrrrrrrrdd]&&&&&&&&\\
\\
&&&&&&&&-1\ar@{-}[lllldd]\ar@{-}[rrrrdd]
&&&&&&&&&&&&&&&&1\ar@{-}[lllldd]\ar@{-}[rrrrdd]&&&&\\
&&&&&&&&&&&&&&&&&&&&&&&&\\
&&&&-2\ar@{-}[lldd]\ar@{-}[rrdd]
&&&&&&&&-\frac{1}{2}\ar@{-}[lldd]\ar@{-}[rrdd]
&&&&&&&&\frac{1}{2}\ar@{-}[lldd]\ar@{-}[rrdd]
&&&&&&&&2\ar@{-}[lldd]\ar@{-}[rrdd]\\
&&&&&&&&&&&&&&&&&&&&&&&&&&&&\\
&&-3\ar@{-}[ldd]\ar@{-}[rdd]
&&&&-\frac{3}{2}\ar@{-}[ldd]\ar@{-}[rdd]
&&&&-\frac{3}{4}\ar@{-}[ldd]\ar@{-}[rdd]
&&&&-\frac{1}{4}\ar@{-}[ldd]\ar@{-}[rdd]
&&&&\frac{1}{4}\ar@{-}[ldd]\ar@{-}[rdd]
&&&&\frac{3}{4}\ar@{-}[ldd]\ar@{-}[rdd]
&&&&\frac{3}{2}\ar@{-}[ldd]\ar@{-}[rdd]
&&&&3\ar@{-}[ldd]\ar@{-}[rdd]\\
&&&&&&&&&&&&&&&&&&&&&&
&&&&&&&&\\
&&
&&&&
&&&&
&&&&
&&&&
&&&&
&&&&&&&&&&\\
&&&&\ldots&&&&&&&&&&&&\ldots&&&&&&&&&&&&\ldots
}
$$
The completion of this picture, according to Conway, contains not only all real numbers, but much more. Conway called surreal numbers an ``enormous
new world of numbers'' and valued them as his greatest discovery; see~\cite[p. 568]{Sch}. Once again, though, we do not know what $q$-deformed
surreal numbers are, but it is tempting to apply the same quantization procedure to the binary tree.  We are not aware, however, of any stabilization
phenomenon.

At this point, John would understand that any sophisticated question he might ask us would result in the same answer \dots . He would probably look
at us with indulgence and ask something simple like: ``But how about quantum complex numbers?! What is the $q$-deformation of~$i$, or
of~$\frac{1+i\sqrt{3}}{2}$?!'' ``Yes," we would reply, ``we {\it do} know the answer to that!:
$$\left[i\right]_q=\frac{i}{\sqrt{q}},$$ while $\frac{1+i\sqrt{3}}{2}$ remains undeformed,
like~$0$ and~$1$," and we would add, ``we'll will tell you more about it very very soon!''

\bigskip
\noindent {\bf Acknowledgements}. We are grateful to Alexander Veselov for enlightening discussions and collaboration. We are pleased to thank Karin
Baur, Justin Lasker, Ludivine Leclere, Bruce Sagan, Richard Schwartz, Michael Shapiro, and Sergei Tabachnikov for fruitful discussions. This paper
was partially supported by the ANR project ANR-19-CE40-0021.


\begin{thebibliography}{99}
\bibitem{Bau}
K.~Baur
{\it Frieze patterns of integers}, arXiv:2101.05676.

\bibitem{Con}
J.~H.~Conway
{\it An enumeration of knots and links, and some of their algebraic properties},
1970 Computational Problems in Abstract Algebra (Proc. Conf., Oxford, 1967), 329--358.

\bibitem{CoCo}
J.~H.~Conway, H.~S.~M.~Coxeter,
{\it Triangulated polygons and frieze patterns,\/}
Math.\ Gaz.\ {\bf 57} (1973), 87--94 and 175--183.

\bibitem{CG}
J.~H.~Conway, R.~K.~Guy, The book of numbers. Copernicus, New York, 1996.

\bibitem{CoR}
H. S. M. Coxeter et J. F. Rigby,
{\it Frieze Patterns, Triangulated Polygons and Dichromatic Symmetry,\/}
The Lighter Side of Mathematics, pp. 15--27, 1994.

\bibitem{Cox}
H.~S.~M.~Coxeter.
{\it Frieze patterns},  Acta Arith.\ {\bf 18}  (1971), 297--310.

\bibitem{CF}
V. Fock, L. Chekhov,
{\it Quantum Teichm\"uller spaces},
Theoret. and Math. Phys. {\bf 120} (1999), 1245--1259.

\bibitem{Kog}
T. Kogiso,
{\it $q$-Deformations and $t$-deformations of Markov triples},
arXiv:2008.12913.

\bibitem{LMG}
L. Leclere, S. Morier-Genoud,
{\it $q$-deformations of the modular group and of the real quadratic irrational numbers},
arXiv:2101.02953.

\bibitem{Ves}
L. Leclere, S. Morier-Genoud, V.~Ovsienko, A. Veselov,
{\it On radius of convergence of q-deformed real numbers,}
arXiv:2102.00891.

\bibitem{R2}
K. Lee, R. Schiffler,
{\it Cluster algebras and Jones polynomials},
Selecta Math. (N.S.) {\bf 25} (2019), no. 4, Art. 58, 41 pp.

\bibitem{CSS}
T. McConville, B. E. Sagan, C. Smyth,
{\it On a rank-unimodality conjecture of Morier-Genoud and Ovsienko},
arXiv:2008.13232.

\bibitem{Sop}
S.~Morier-Genoud,
{\it Coxeter's frieze patterns at the crossroads of algebra, geometry and combinatorics},
Bull.\ Lond.\ Math.\ Soc.\ {\bf 47} (2015), 895--938.

\bibitem{MGO}
S.~Morier-Genoud, V.~Ovsienko,
{\it Farey boat: continued fractions and triangulations, modular group and polygon dissections,\/}
Jahresber. Dtsch. Math.-Ver. {\bf 121} (2019), no. 2, 91--136.

\bibitem{SVRat}
S. Morier-Genoud, V. Ovsienko,
{\it $q$-deformed rationals and $q$-continued fractions.}
Forum Math. Sigma 8 (2020), e13, 55 pp.

\bibitem{SVRe}
S.~Morier-Genoud, V.~Ovsienko,
{\it On q-Deformed Real Numbers},
Experimental Mathematics, DOI: 10.1080/10586458.2019.1671922, arXiv:1908.04365.

\bibitem{OEIS}
OEIS Foundation Inc., The On-Line Encyclopedia of Integer Sequences, http://oeis.org.

\bibitem{Sch}
D. Schleicher,
{\it Interview with John Horton Conway}, Notices Amer. Math. Soc.
{\bf 60} (2013), no. 5, 567--575.

\end{thebibliography}
\end{document}